\documentclass[12pt]{article}
\usepackage{color}
\definecolor{darkblue}{rgb}{0.00,0.05,0.50}
\usepackage[colorlinks,filecolor=blue,citecolor=darkblue]{hyperref}

\usepackage{xcolor}
\usepackage{hyperref}

\definecolor{linkcolor}{HTML}{00FF00} % колір посилань
\definecolor{urlcolor}{HTML}{00FF00} % колір гіперпосилань

\hypersetup{pdfstartview=FitH,  linkcolor=linkcolor,urlcolor=urlcolor, colorlinks=true}

\usepackage{amsfonts,amssymb,amsmath,amsthm}
\usepackage{mathrsfs}
\usepackage{url}
\usepackage{enumerate}
\usepackage[ukrainian, russian, english]{babel}
\usepackage[cp1251]{inputenc}
\begin{document} \selectlanguage{ukrainian}
\thispagestyle{empty}

\title{}

UDC 517.51 \vskip 5mm

\begin{center}
\textbf{\Large Estimates for the entropy numbers of the Nikol'skii--Besov classes of periodic functions of many variables in the space of quasi-continuous functions}
\end{center}

\vskip 3mm

\begin{center}
\textbf{\Large   Оцінки ентропійних чисел класів Нікольського--Бєсова періодичних функцій багатьох змінних у просторі квазінеперервних функцій}
\end{center}
\vskip0.5cm

\begin{center}
 A.\,S.~Romanyuk, S.\,Ya.~Yanchenko\\ \emph{\small
Institute of Mathematics of NAS of
Ukraine, Kyiv}
\end{center}
\begin{center}
А.\,С.~Романюк, C.\,Я.~Янченко \\
\emph{\small Інститут математики НАН України, Київ}
\end{center}
\vskip0.5cm

\begin{abstract}

We obtained exact-order estimates for the entropy numbers of the Nikol'skii--Besov classes $B^{\boldsymbol{r}}_{p,\theta}$ of periodic functions of many variables in the metric of the space of quasi-continuous functions.

\vskip 5 mm

Отримано порядкові оцінки енторпійних чисел класів Нікольського--Бєсова $B^{\boldsymbol{r}}_{p,\theta}$ періодичних функцій багатьох змінних у метриці простору квазінеперервних функцій.

\end{abstract}

\vskip 0.7 cm

%%%%%% Вступ %%%%%%%%%%%%%%%%%%%%%%%%%%%%%%%%%%%%%%%%%%%%%%%%%%%%%%%%%%%%%%%%%%%%%%%%

\vskip 1mm  \textbf{1. Вступ}. У роботі продовжено вивчення асимптотичних характеристик класів Нікольського--Бєсова $B^{\boldsymbol{r}}_{p,\theta}$ з домінуючою мішаною похідною періодичних функцій багатьох змінних \cite{Romanyuk_2019AM}, \cite{Romanyuk_2017UMJ}, \cite{Romanyuk_2015UMJ}. Знайдено оцінки для ентропійних чисел згаданих класів у метриці $QC$-простору квазінеперервних функцій. Показано, що у деяких випадках ці оцінки є точними за порядком. Даний простір функцій є близьким за своїми властивостями до простору $L_{\infty}$, а знаходження оцінок такої характеристики, як ентропійні числа, в його метриці дає змогу одержати нові результати~--- досі не встановлені у рівномірній метриці.

\vskip 3mm  \textbf{2. Означення класів функцій та апроксимативних характеристик}. Нехай $\mathbb{R}^d$, $d\geq 1$,~---  евклідів  простір з  елементами $\boldsymbol{x}=(x_1,
 \ldots, x_d)$ і $(\boldsymbol{x},\boldsymbol{y}) = x_1 y_1+ \ldots + x_d y_d$. Через $L_p(\pi_d)$, $\pi_d = \prod\limits^{d}_{j=1} [0,2\pi]$, $1\leq p \leq \infty$, позначимо простір функцій $f(\boldsymbol{x})$, які є  $2\pi$-періодичні за кожною змінною зі скінченною нормою:
  $$
  \|f\|_p : = \|f\|_{L_p(\pi_d)} =  \left((2\pi)^{-d}
  \int\limits_{\pi_d}|f(\boldsymbol{x})|^p\,d\boldsymbol{x}\right)^{1/p} <\infty, 1 \leq p < \infty,
 $$
 $$
 \|f\|_{\infty}:=\|f\|_{L_{\infty}(\pi_d)}=\mathop {\rm ess \sup}\limits_{\boldsymbol{x}\in \pi_d} |f(\boldsymbol{x})|.
 $$

У подальших міркуваннях будемо розглядати тільки ті функції  ${f \in L_p(\pi_d)}$, для яких виконана умова
$$
\int\limits^{2\pi}_0 f(\boldsymbol{x})dx_j = 0,  j=\overline{1,d},
$$
і множину таких функцій будемо позначати $L^0_p(\pi_d)$.

Для функції $f \in L^0_p(\pi_d)$, $1 \leq p \leq \infty$,   розглянемо різницю першого порядку  по $j$-ій  змінній з кроком $h\in \mathbb{R}$:
$$
 \Delta_{h,j}f(\boldsymbol{x})=f(x_1,\dots,x_{j-1},x_j+h,x_{j+1},\dots,x_d)-f(\boldsymbol{x})
$$
і, відповідно, $l$-го порядку, $l \in \mathbb{N}$,
$$
\Delta_{h,j}^{l}f(\boldsymbol{x})=
 \overbrace{\Delta_{h,j}\dots\Delta_{h,j}}\limits^{l}f(\boldsymbol{x}).
$$

Далі, якщо $\boldsymbol{k} = (k_1, \ldots , k_d)$, $k_j \in \mathbb{N}$, $j = \overline{1,d}$,  то  мішана різниця порядку  $\boldsymbol{k}$ з векторним кроком $\boldsymbol{h} = (h_1, \ldots, h_d)$, $h_j\in \mathbb{R}$, $j = \overline{1,d}$, визначається рівністю
$$
 \Delta_{\boldsymbol{h}}^{\boldsymbol{k}}f(\boldsymbol{x})=
 \Delta_{h_1,1}^{k_1}\Delta_{h_2,2}^{k_2}\dots\Delta_{h_d,d}^{k_d}f(\boldsymbol{x}).
$$

Нехай задані  вектор $\boldsymbol{r} = (r_1, \ldots, r_d)$, $r_j > 0$ $j = \overline{1,d}$,  і параметри $1 \leq p, \theta \leq \infty$. Тоді простори  $B^{\boldsymbol{r}}_{p,\theta}(\pi_d)$ можна означити таким чином:
$$
 B^{\boldsymbol{r}}_{p,\theta}:=B^{\boldsymbol{r}}_{p,\theta}(\pi_d)=\Big \{f \in
 L^0_p(\pi^d)\colon \|f\|_{B^{\boldsymbol{r}}_{p,\theta}} <\infty \Big\},
$$
і норма задається рівностями
$$
  \|f\|_{B^{\boldsymbol{r}}_{p,\theta}} = \left( \int\limits_{\pi_d}\|\Delta^{\boldsymbol{k}}_{\boldsymbol{h}} f\|^{\theta}_p  \prod\limits^d_{j=1}\frac{d h_j}{h^{1+r_j \theta}_j} \right)^{1/\theta}
$$                                                                                                                                                                                    якщо $1\leqslant \theta < \infty$ та
$$
\|f\|_{H^{\boldsymbol{r}}_{p}}\equiv\|f\|_{B^{\boldsymbol{r}}_{p,\infty}} = \sup\limits_{\boldsymbol{h}} \|\Delta^{\boldsymbol{k}}_{\boldsymbol{h}} f\|_p \prod\limits^d_{j=1}h^{-r_j}_j.
$$
Також вважаємо, що  для векторів  $\boldsymbol{k} = (k_1, \ldots, k_d)$  і $\boldsymbol{r} = (r_1, \ldots, r_d)$ виконана умова $k_j > r_j$, $j = \overline{1,d}$.

У такій формі означення просторів $B^{\boldsymbol{r}}_{p,\theta}$ було дане у роботах В.\,М.~Темлякова~\cite{Temlyakov86} і С.\,М.~Нікольського та П.\,І.~Лізоркіна~\cite{Lizorkin_Nikolsky_1989} відповідно для $H^{\boldsymbol{r}}_p$ і $B^{\boldsymbol{r}}_{p,\theta}$. Вони належать шкалі просторів мішаної гладкості, що введені С.\,М.~Нікольським~\cite{Nikolsky_63} і Т.\,І.~Амановим~\cite{Amanov_1965}. Окрім того, вони є узагальненням відомих ізотропних просторів Бєсова~\cite{Besov_1961}, а для випадку $\theta=\infty$~--- Нікольського~\cite{Nikolsky_51}.

Під класом $B^{\boldsymbol{r}}_{p,\theta}$ будемо розуміти множину
функцій ${f \in L^0_p(\pi^d)}$ для яких
${\|f\|_{B^{\boldsymbol{r}}_{p,\theta}}\leqslant  1}$ і при цьому збережемо
для класів $B^{\boldsymbol{r}}_{p,\theta}$ ті ж самі позначення, що і для
просторів $B^{\boldsymbol{r}}_{p,\theta}$.

При проведенні подальших міркувань  нам буде зручно користуватися  означенням норми функцій з  класів $B^{\boldsymbol{r}}_{p,\theta}$ в дещо іншій формі, а саме опосередковано через, так зване, декомпозиційне представлення елементів цих просторів. Уперше декомпозиційне представлення функцій з класів Нікольського--Бєсова та відповідне йому нормування з'явилося у роботах В.\,М.~Темлякова~\cite[Роз.~2, п.\,1]{Temlyakov86} і  С.\,М.~Нікольського та П.\,І.~Лізоркіна~\cite{Lizorkin_Nikolsky_1989} і, як з'ясувалося пізніше, відіграло ключову роль у дослідженнях, які пов'язані з апроксимацією класів функцій.

Для векторів
$\boldsymbol{s} = (s_1, \ldots, s_d)$,  $s_j \in \mathbb{Z}_+$,  $\boldsymbol{k}= (k_1, \ldots, k_d)$, $k_j \in \mathbb{Z}$, $j=\overline{1,d}$, покладемо
$$
\rho(\boldsymbol{s}) = \big\{ \boldsymbol{k} = (k_1, \ldots, k_d)\colon [2^{s_j - 1}] \leq |k_j| < 2^{s_j}, j=\overline{1,d} \big\}
$$
 і для $f \in L^0_p(\pi_d)$  позначимо
 $$
 \delta_{\boldsymbol{s}}(f):= \delta_{\boldsymbol{s}}(f, \boldsymbol{x}) = \sum\limits_{\boldsymbol{k} \in \rho(\boldsymbol{s})} \widehat{f}(\boldsymbol{k}) e^{i(\boldsymbol{k},\boldsymbol{x})},
 $$
де
$\widehat{f}(\boldsymbol{k}) =  \int\limits_{\pi_d} f(\boldsymbol{t}) e^{-i(\boldsymbol{k},\boldsymbol{t})}d\boldsymbol{t}$~--- коефіцієнти Фур'є функції $f$.

Тоді класи $B_{p,\theta}^{\boldsymbol{r}}$,
$1<p<\infty$, $1\leqslant  \theta \leqslant \infty$, $\boldsymbol{r}= (r_1, \ldots, r_d)$ $r_j>0$ $j = \overline{1,d}$, можна означити наступним
чином~\cite{Lizorkin_Nikolsky_1989}:
$$
B^{\boldsymbol{r}}_{p,\theta}:=\Big\{f\in L^0_p(\pi_d)\colon \
\|f\|_{B^{\boldsymbol{r}}_{p,\theta}} \leqslant 1 \Big\},
$$
де
$$
   \|f\|_{B^{\boldsymbol{r}}_{p,\theta}}\asymp \Bigg(\sum \limits_{\boldsymbol{s}\in\mathbb{Z}^d_+}2^{(\boldsymbol{s},\boldsymbol{r})\theta}
   \|\delta_{\boldsymbol{s}}(f)\|_p^{\theta}\Bigg)^{\frac{1}{\theta}}
$$
  при $1\leqslant\theta<\infty$ і
$$
   \|f\|_{B^{\boldsymbol{r}}_{p,\infty}}\equiv\|f\|_{H^{\boldsymbol{r}}_{p}}\asymp \sup
   \limits_{\boldsymbol{s}\in\mathbb{Z}^d_+} 2^{(\boldsymbol{s},\boldsymbol{r})}\|\delta_{\boldsymbol{s}}(f)\|_p.
$$

Тут і надалі по тексту для додатних величин $a$ і  $b$ вживається запис
$a\asymp b$, який означає, що існують такі додатні  сталі $C_1$ та $C_2$,
які не залежать від одного істотного параметра у величинах  $a$ і  $b$,
що $C_1 a \leqslant b$ (у цьому випадку пишемо $a\ll b$) і  $C_2 a \geqslant b$ (у цьому випадку пишемо $a\gg b$).  Всі сталі $C_i$, $i=1,2,\dots$, які зустрічаються у роботі, можуть залежати  лише від тих параметрів, що
входять в означення класу, метрики, в якій оцінюється похибка
наближення, та розмірності простору $\mathbb{R}^d$. У деяких випадках  цю  залежність  будемо вказувати у явному вигляді.

Зазначимо, що видозмінивши ``блоки''  $\delta_{\boldsymbol{s}}(f)$, наведене означення класів $B^{\boldsymbol{r}}_{p,\theta}$  можна поширити і на  крайні значення $p=1$ і $p=\infty$ (див., наприклад, \cite[зауваження~2.1]{Lizorkin_Nikolsky_1989}.

Нехай $V_l(t)$,  $t\in\mathbb{R}$, $l \in \mathbb{N}$, позначає ядро Валле-Пуссена вигляду
$$
V_l(t) = 1 + 2 \sum\limits^l_{k=1} \cos kt  + 2\sum\limits^{2l-1}_{k=l+1}\left(1 - \frac{k-l}{l}\right)\cos kt.
$$
Поставимо у відповідність кожному вектору $\boldsymbol{s} = (s_1, \ldots, s_d)$, $s_j \in \mathbb{Z}_+$, $j=\overline{1,d}$, поліном
$$
A_{\boldsymbol{s}}(\boldsymbol{x}) = \prod\limits^d_{j=1}(V_{2^{s_j}}(x_j) - V_{2^{s_j - 1}}(x_j))
$$
 і для $f \in L^0_p(\pi_d)$, $1 \leq p \leq \infty$,  покладемо
 $$
 A_{\boldsymbol{s}}(f) := A_{\boldsymbol{s}}(f,\boldsymbol{x}) = (f \ast A_{\boldsymbol{s}})(\boldsymbol{x}),
 $$
де  ``$\ast$'' означає операцію згортки.   Тоді, при $1 \leq p \leq \infty$,  $\boldsymbol{r} = (r_1, \ldots, r_d)$, $r_j > 0$, $j = \overline{1,d}$, класи $B^{\boldsymbol{r}}_{p,\theta}$  можна означити таким чином:
$$
B^{\boldsymbol{r}}_{p,\theta} = \Big\{ f\colon \|f\|_{B^{\boldsymbol{r}}_{p,\theta}} \asymp \Bigg(\sum\limits_{\boldsymbol{s}\in \mathbb{Z}^d_+} 2^{(\boldsymbol{s},\boldsymbol{r})\theta} \|A_{\boldsymbol{s}}(f)\|^{\theta}_p \Bigg)^{\frac{1}{\theta}} \leq 1, 1 \leq \theta < \infty  \Big\},
$$
$$
B^{\boldsymbol{r}}_{p,\infty} = \Big\{ f\colon \|f\|_{B^{\boldsymbol{r}}_{p,\infty}} \asymp \sup\limits_{s\in \mathbb{Z}^d_+} 2^{(\boldsymbol{s},\boldsymbol{r})} \|A_{\boldsymbol{s}}(f)\|_p \leq 1 \Big \}.
$$

З дослідженням різних апроксимативних характеристик класів Нікольського та Нікольського--Бєсова періодичних функцій можна ознайомитися у монографіях \cite{Cross_2018}, \cite{Romanjuk 2012}, \cite{Temlyakov86}, \cite{Temlyakov 1993},  де наведена детальна бібліографія.

Тепер означимо  асимптотичні характеристики, які будемо досліджувати.

Нехай $\mathscr{X}$  банахів простір і $B_{\mathscr{X}}(\boldsymbol{y},r) = \big\{ x \in \mathscr{X}\colon \|\boldsymbol{x} - \boldsymbol{y} \| \leq r \big\}$~--- куля радіуса $r$ з центром у точці $\boldsymbol{y}$.

Для компактної  множини $A\subset \mathscr{X}$  і $\varepsilon > 0$  позначимо
$$
N_{\varepsilon}(A,\mathscr{X}) = \min \Big\{n \colon \exists \boldsymbol{y}^1, \ldots, \boldsymbol{y}^n \in \mathscr{X} \colon A\subseteq \bigcup\limits^n_{j=1} B_{\mathscr{X}}(\boldsymbol{y}^j, \varepsilon)\Big\}.
$$
Тоді величина
$$
H_{\varepsilon}(A,\mathscr{X}) = \log N_{\varepsilon}(A,\mathscr{X})
$$
називається  $\varepsilon$-ентропією множини $A$ відносно  банахового  простору $\mathscr{X}$~\cite{Kolmogorov_Tixomirov-59} (тут і далі під записом $\log$ будемо розуміти $\log_2$).

З $\varepsilon$-ентропією множини  $A$ тісно пов'язано  поняття  її  ентропійних  чисел $\varepsilon_k(A,\mathscr{X})$ (див., наприклад, \cite{Hollig 1980}):
$$
\varepsilon_k(A,\mathscr{X}) = \inf \Big\{\varepsilon \colon \exists \boldsymbol{y}^1, \ldots, \boldsymbol{y}^{2^k} \in \mathscr{X} \colon A\subseteq \bigcup\limits^{2^k}_{j=1} B_{\mathscr{X}}(\boldsymbol{y}^j,\varepsilon) \Big \}.
$$
Безпосередньо з означень величин $H_{\varepsilon}(A,\mathscr{X})$ і $\varepsilon_k(A,\mathscr{X})$ можемо записати: якщо $H_{\varepsilon}(A,\mathscr{X}) \leq k$,  то $\varepsilon_k(A,\mathscr{X}) \leq \varepsilon$ і навпаки~--- з оцінки  ${{\varepsilon}_k(A,\mathscr{X}) \leq \varepsilon}$ отримуємо ${H_{\varepsilon}(A,\mathscr{X}) \leq k}$, а саме, якщо $k< H_{\varepsilon}(A,\mathscr{X}) \leq k+1$, то $\varepsilon_{k+1}(A,\mathscr{X})\leq \varepsilon \leq \varepsilon_k(A,\mathscr{X})$. Ці  співвідношення дають  можливість  із оцінок  для ентропійних чисел деякої множини $A$ отримувати відповідні оцінки її $\varepsilon$-ентропії.

Дослідження $\varepsilon$-ентропії і близьких до неї асимптотичних характеристик ($\varepsilon$-ємність, ентропійні числа і т.п.) мають багату історію. Зокрема, ентропійні числа для класів функцій однієї та багатьох змінних Соболєва $W^{\boldsymbol{r}}_{p,\alpha}$, Нікольського $H^{\boldsymbol{r}}_{p}$, Нікольського--Бєсова $B^{\boldsymbol{r}}_{p,\theta}$ та їх аналогів досліджувалися у роботах \cite{Belinskii_1989AM}, \cite{Belinskii_1990Ya},  \cite{Belinskii_1998JAT}, \cite{Dung_2001JCom},  \cite{Dunker_1999JAT},  \cite{Kashin_Temlyakov_1995MN},  \cite{Kashin_Temlyakov_1998MN}, \cite{Kashin_Temlyakov_1999},  \cite{Kashin_Temlyakov_1994MN},  \cite{Mayer_Ullrich-2020},  \cite{Pozharska_2019JMSci}, \cite{Pozharska_2019UMJ}, \cite{Romanyuk_2019AM}, \cite{Romanyuk_2017UMJ}, \cite{Romanyuk_2015UMJ}, \cite{RomanyukASVS_2019UMJ}, \cite{Temlyakov_2013JAT}, \cite{Temlyakov_1995JC},  \cite{Temlyakov_96EJA}, \cite{Temlyakov_1989TR}, \cite{Temlyakov_2017JAT}, \cite{Temlyakov_1998EJA},  \cite{Temlyakov_Ullrich-arx2020},  \cite{Vybiral_2006}. З більш детальною  бібліографією можна ознайомитися у монографіях \cite{Cross_2018}, \cite{Temlyakov_2018}, \cite{Trigub_Belinsky}.

Означимо простір у метриці якого будемо оцінювати ентропійні числа.

Нехай $\mu$~--- нормована міра Лебега на одиничному колі. Для функції $f \in L_1(d\mu)$ з рядом Фур'є
$$
f\sim \sum\limits_{s=0}^{\infty}\delta_{s}(f, x),
$$
$$
\delta_{0}(f, x)=\int^{2\pi}_0 f d\mu, \ \ \delta_{s}(f, x) = \sum\limits_{2^{s-1}\leq |k| < 2^s} \widehat{f}(k) e^{ikx}, \ s=1, 2, \dots  ,
$$
розглянемо величину
\begin{equation}\label{QC-norm}
\|f\|_{QC}\equiv \int^1_0 \bigg\| \sum^{\infty}_{s=0} r_s(\omega)\delta_{s}(f, x)\bigg\|_{L_{\infty}(d\mu)}d\omega,
\end{equation}
де $\big\{ r_s(\omega)\big\}^{\infty}_{s=0}$~--- система Радемахера (див.,~\cite[Гл.\,2, \S\,1]{Kashin-Saakyan_1984}).  Тоді простором квазінеперервних функцій (позначення $QC$) будемо називати замикання множини тригонометричних поліномів за нормою \eqref{QC-norm}.

Простір квазінеперервних функцій у багатовимірному випадку ($d\geq 2$) означимо таким чином:
\begin{equation}\label{QC-dnorm}
\|f\|_{QC}\equiv \big\|\|f(\cdot, \boldsymbol{x}^1)\|_{QC}\big\|_{\infty},
\end{equation}
де для $\boldsymbol{x}=(x_1,\ldots,x_d)\in \pi_d$ покладаємо $\boldsymbol{x}^1=(x_2,\ldots,x_d)\in \pi_{d-1}$, тобто в \eqref{QC-dnorm} береться $QC$-норма по змінній $x_1$ і $\sup$-норма по інших змінних.

Зазначимо, що простір квазінеперервних функцій $QC$ введено в роботі \cite{Kashin_Temlyakov_1998MN} (див. також \cite{Kashin_Temlyakov_1999}), де наведено деякі його властивості, зокрема, при $d=1$ для $f\in QC$ справедливе співвідношення
$$
\|f\|_{QC}\leq \sum\limits_s \|\delta_s(f)\|_{\infty}.
$$
З дослідженнями властивостей, а також із застосуваннями $QC$-норми можна ознайомитися у роботах \cite{Radomskii2016-MS}, \cite{Radomskii2018-UMN}.

Якщо    $\mathfrak{M}$ --- деяка скінченна множина, то через $|\mathfrak{M}|$  будемо позначати кількість її елементів.

\vskip 3mm

\textbf{3. Основні результати.} У подальших міркування будемо вважати, що вектор $\boldsymbol{r}$, який входить в означення класів $B^{\boldsymbol{r}}_{p,\theta}$ має вигляд $\boldsymbol{r}=(r_1, \ldots, r_1)\in \mathbb{R}_+^d$.

\vskip 1 mm
\bf Теорема 1.  \rm  \it Нехай $1 < p \leq \infty$,  $1 \leq \theta < \infty$, $r_1 >\max\left\{\frac{1}{p}, \frac{1}{2}\right\}$. Тоді при $d\geq 2$ справедлива оцінка
\begin{equation}\label{Em-QC}
\varepsilon_M\big(B^{\boldsymbol{r}}_{p,\theta}, QC\big) \ll  M^{-{r_1}}(\log^{d - 1} M)^{r_1 +( \frac{1}{p^*} - \frac{1}{\theta})_{+}} \sqrt{\log\, M},
\end{equation}
де $p^*=\min \{p, 2\}$, $a_{+} = \max \{ a,  0\}$. \rm

{\textbf{\textit{Доведення.}}} Розглянемо спочатку випадок $1<p \leq 2$, $p<\theta<\infty$.

Для $n \in \mathbb{N}$ покладемо
$$
Q_n = \bigcup\limits_{(\boldsymbol{s},1)\leq n} \rho(\boldsymbol{s}), \Delta Q_n =  Q_n \backslash Q_{n-1} \ \ \text{і} \ \ \mathfrak{N}_n = \big \{ \boldsymbol{s} = (s_1, \ldots, s_d), (\boldsymbol{s}, 1) = n\big\}.
$$
Зазначимо, що $|\Delta Q_n | \asymp 2^n n^{d - 1}$.

Для подальших міркувань нам знадобляться такі допоміжні твердження.

\textbf{Лема А} (див., наприклад, \cite{Temlyakov 1993})\textbf{.} \it Нехай $ f \in L^0_p(\pi_d)$,  $1 < p < \infty$.  Тоді
$$
\bigg\|\sum\limits_{\boldsymbol{s}} \delta_{\boldsymbol{s}}(f)\bigg\|_p \ll \left( \sum\limits_{\boldsymbol{s}} \|\delta_{\boldsymbol{s}}(f)\|^{p^*}_p\right)^{\frac{1}{p^*}},
$$
де $p^* = \min \{p, 2\}$. \rm

\textbf{Лема Б} \cite[с.\,11]{Temlyakov86}\textbf{.} \it Справедливе співвідношення
$$
\sum\limits_{(\boldsymbol{s}, 1) \geq l}  2^{-{\alpha}(\boldsymbol{s},1)} \asymp 2^{-{\alpha} l} l^{d - 1},  \ \ \alpha > 0.
$$\rm

Отже, згідно з лемою~А для $f\in B^{\boldsymbol{r}}_{p,\theta}$ будемо мати
\begin{equation}\label{delta-es1}
 \bigg \| \sum \limits_{\boldsymbol{s} \in \mathfrak{N}_n } \delta_{\boldsymbol{s}}(f) \bigg \|_p \ll
 \left( \sum \limits_{\boldsymbol{s} \in \mathfrak{N}_n}\|\delta_{\boldsymbol{s}} (f)\|_{p}^p \right)^{\frac
{1}{p}} = \mathfrak{I}_1.
 \end{equation}

Далі, скориставшись нерівністю Гельдера з показником $\theta/p$ і лемою~Б, одержимо
$$
\mathfrak{I}_1 \leq \left( \sum \limits_{\boldsymbol{s} \in \mathfrak{N}_n } 2^{(\boldsymbol{s},\boldsymbol{r})}\|\delta_{\boldsymbol{s}}(f)\|_{p}^{\theta}\right)^{\frac
{1}{\theta}} \left( \sum \limits_{\boldsymbol{s} \in \mathfrak{N}_n } 2^{-(\boldsymbol{s},\boldsymbol{r})\frac
{\theta p}{\theta -p}}\right)^{\frac
{1}{p}-\frac{1}{\theta}}\ll
$$
$$
\ll \|f\|_{B_{p,\theta}^{\boldsymbol{r}}} \left( \sum \limits_{\boldsymbol{s} \in \mathfrak{N}_n } 2^{-(\boldsymbol{s},\boldsymbol{r})\frac
{\theta p}{\theta -p}}\right)^{\frac
{1}{p}-\frac{1}{\theta}} \leq
$$
\begin{equation}\label{delta-es11}
\leq \left( \sum \limits_{\boldsymbol{s} \in \mathfrak{N}_n } 2^{-r_1(\boldsymbol{s},1)\frac
{\theta p}{\theta -p}}\right)^{\frac
{1}{p}-\frac{1}{\theta}} \asymp 2^{-n r_1} n^{(d-1)\left(\frac
{1}{p}-\frac{1}{\theta}\right)}.
 \end{equation}

Аналогічно у випадку $1<p\leq 2$ і $\theta=p$ можемо записати
$$
\mathfrak{I}_1 \ll 2^{-n r_1} \left( \sum \limits_{\boldsymbol{s} \in \mathfrak{N}_n } 2^{(\boldsymbol{s},\boldsymbol{r})}\|\delta_{\boldsymbol{s}}(f)\|_{p}^{p}\right)^{\frac
{1}{p}}\ll
$$
\begin{equation}\label{delta-theta=p}
\ll 2^{-n r_1} \|f\|_{B_{p,p}^{\boldsymbol{r}}} \ll 2^{-n r_1}.
\end{equation}
Таким чином, при $1<p\leq2$, $p\leq \theta < \infty$, згідно з \eqref{delta-es1}--\eqref{delta-theta=p} маємо
\begin{equation}\label{delta-fin}
\bigg \| \sum \limits_{\boldsymbol{s} \in \mathfrak{N}_n } \delta_{\boldsymbol{s}}(f) \bigg \|_p \ll 2^{-n r_1} n^{(d-1)\left(\frac
{1}{p}-\frac{1}{\theta}\right)}.
\end{equation}

Далі, нехай задано достатньо велике число $M$. Підберемо $m\in \mathbb{N}$ так, щоб виконувалися нерівності
$$
|\Delta Q_{m-1}| < M \leq |\Delta Q_{m}|.
$$

Тоді, оскільки
$$
|\Delta Q_{m-1}| \asymp |\Delta Q_{m}| \asymp 2^m m^{d-1},
$$
то $M\asymp 2^m m^{d-1}$.

Тепер покладемо $\beta =  \frac{1}{2}\min \left\{\left(r_1 - \frac{1}{p}\right), 1 \right\}$  і
$$
\overline{M}_n = \left\{
\begin{array}{ll}
   C(\beta) M 2^{-\frac{1}{2}(m - n)}, n < m,  \\
   C(\beta) M 2^{-{\beta}(n - m)}, n \geq m,
\end{array} \right.
$$
де числа $C(\beta) > 0$ підібрано таким чином, що
$$
\sum\limits^{\infty}_{n=1} \overline{M}_n \leq M.
$$
Зауважимо, що такі числа $C(\beta) > 0$ існують, оскільки
$$
M \sum\limits_{n=0}^{m-1} 2^{-\frac{1}{2}(m-n)} + M \sum\limits_{n=m}^{\infty}   2^{-{\beta}(n-m)}  \ll M.
$$

Нехай $M_n = [\overline{M}_n]$,  де $[a]$~---  ціла частина числа $a$. Тоді ${M_n = 0}$, якщо  $C(\beta) M 2^{-{\beta}(n-m)} < 1$, тобто при ${n> m_1 = m + {\beta}^{-1}\log C(\beta) M}$.

Позначимо
$$
S_{\Delta Q_n} \big(B^{\boldsymbol{r}}_{p,\theta}\big) := \Big \{ g \colon g(\boldsymbol{x}) = \sum\limits_{\boldsymbol{k} \in \Delta Q_n} \widehat{f}(\boldsymbol{k}) e^{i(\boldsymbol{k},\boldsymbol{x})}, f \in B^{\boldsymbol{r}}_{p,\theta} \Big\}
$$
і
$$
\big\|S_{\Delta Q_n}\big(B^{\boldsymbol{r}}_{p,\theta}\big) \big \|_{QC} : = \sup\limits_{g \in S_{\Delta Q_n}(B^{\boldsymbol{r}}_{p,\theta})} \|g\|_{QC}.
$$

Отже, згідно з позначеннями для ентропійних чисел  $\varepsilon_M\big(B^{\boldsymbol{r}}_{p,\theta}, QC \big)$, можемо записати
\begin{equation}\label{epsilon-sum}
\varepsilon_M\big(B^{\boldsymbol{r}}_{p,\theta}, QC \big) \leq \sum\limits_{n \leq m_1} \varepsilon_{M_n}\big(S_{\Delta Q_n}(B^{\boldsymbol{r}}_{p,\theta}), QC \big) +
\sum\limits_{n > m_1} \| S_{\Delta Q_n}(B^r_{p,\theta})\|_{QC}  = I_1 + I_2.
\end{equation}

Оцінимо спочатку доданок $I_2$, скориставшись відомим твердженням.

Для будь-якої множини $\Lambda \subset \mathbb{Z}^d$ через $\mathcal{T}(\Lambda)$ будемо позначати множину тригонометричних поліномів $t$ вигляду
$$
t(\boldsymbol{x}) = \sum\limits_{\boldsymbol{k} \in \Lambda} c_{\boldsymbol{k}} e^{i(\boldsymbol{k},\boldsymbol{x})}, \boldsymbol{x} \in \pi_d.
$$
У випадку, коли множина $\Lambda$ симетрична відносно початку координат $(\Lambda = -  \Lambda)$, покладемо
$$
\mathcal{T}_r(\Lambda) = \big\{t \in \mathcal{T}(\Lambda)\colon c_{\boldsymbol{k}} = \overline{c}_{-\boldsymbol{k}}, \boldsymbol{k} \in \Lambda \big\}.
$$

\bf Теорема А \rm \cite[ Роз.\,1, теорема~2.1]{Temlyakov86}\textbf{.} \it  Нехай $f \in \mathcal{T}(Q_n)$. Тоді при  $1\leq p < \infty$ справедлива нерівність
\begin{equation}\label{TeorA}
  \|f\|_{\infty} \ll 2^{\frac{n}{p}} n^{(d - 1)\left(1 - \frac{1}{p}\right)} \|f\|_{p}.
\end{equation}\rm

\vskip 3 mm

Зауважимо, що оцінка \eqref{TeorA} залишається вірною і в тому випадку, коли ${f\in\mathcal{T}(\Delta Q_n)}$.

Оцінимо спочатку $\|g\|_{QC}$.  Згідно з означенням та властивостями $QC$-норми і оцінками \eqref{delta-fin} та \eqref{TeorA} маємо
$$
\|g\|_{QC}=\big\| \|g(\cdot, \boldsymbol{x}^1)\|_{QC}\big\|_{\infty}\ll
$$
$$
\ll \bigg \| \sum \limits_{\boldsymbol{s} \in \mathfrak{N}_n } \delta_{\boldsymbol{s}}(f) \bigg \|_{\infty} \ll 2^{\frac{n}{p}} n^{(d - 1)\left(1 - \frac{1}{p}\right)}  \bigg \| \sum \limits_{\boldsymbol{s} \in \mathfrak{N}_n } \delta_{\boldsymbol{s}}(f) \bigg \|_{p} \ll
$$
$$
\ll 2^{\frac{n}{p}} n^{(d - 1)\left(1 - \frac{1}{p}\right)} 2^{-n r_1} n^{(d-1)\left(\frac
{1}{p}-\frac{1}{\theta}\right)}  = 2^{-n\left(r_1 - \frac{1}{p}\right)} n^{(d - 1)\left(1 -  \frac{1}{\theta}\right)}.
$$

Отже, для кожного доданку в $I_2$ можемо записати
$$
\big \|S_{\Delta Q_n}\big(B^{\boldsymbol{r}}_{p,\theta}\big) \big \|_{QC} \ll 2^{-n\left(r_1 - \frac{1}{p}\right)} n^{(d - 1)\left(1 -  \frac{1}{\theta}\right)}.
$$

Далі, провівши підсумовування по $n>m_1$ і врахувавши значення $m_1$, одержимо
\begin{equation}\label{I_2os}
I_2  \ll \sum\limits_{n > m_1} 2^{-n\left(r_1 - \frac{1}{p}\right)} n^{(d - 1)\left(1 - \frac{1}{\theta}\right)} \ll  2^{-m_1\left(r_1 - \frac{1}{p}\right)} m_1^{(d - 1)\left(1 -  \frac{1}{\theta}\right)} = \mathfrak{I}_2.
\end{equation}

Для продовження оцінки величини $\mathfrak{I}_2$ розглянемо два випадки.

Нехай $r_1-\frac{1}{p}>1$. Тоді $\beta=\frac{1}{2}$ і відповідно $m_1=m+\log(C(\beta)M)^2$. Таким чином будемо мати
$$
\mathfrak{I}_2 = 2^{-m\left(r_1 -\frac{1}{p}\right)}(C(\beta)M)^{-2\left(r_1 - \frac{1}{p}\right)}\big(m + \log(C(\beta)M)^2\big)^{(d - 1)\left(1 -  \frac{1}{\theta}\right)} \asymp
$$
$$
\asymp 2^{-m\left(r_1 -\frac{1}{p}\right)}2^{-2\left(r_1 -\frac{1}{p}\right)m} m^{-2(d - 1)\left(r_1 -\frac{1}{p}\right)} m^{(d - 1)\left(1 -  \frac{1}{\theta}\right)} \ll
$$
\begin{equation}\label{I_2fin1}
\ll 2^{-m r_1}  m^{(d - 1)\left(\frac{1}{p} - \frac{1}{\theta}\right)}.
\end{equation}

Нехай тепер $0< r_1-\frac{1}{p}\leq 1$. Тоді $\beta=\frac{1}{2}\left(r_1-\frac{1}{p}\right)$,  $m_1=m+\log(C(\beta)M)^{\frac{2p}{r_1 p-1}}$ і величина $\mathfrak{I}_2$ оцінюється у такий спосіб
$$
\mathfrak{I}_2 = 2^{-m\left(r_1 - \frac{1}{p}\right)}(C(\beta)M)^{-\frac{2p}{r_1 p-1}\left(r_1-\frac{1}{p}\right)}\big(m + \log(C(\beta)M)^{\frac{2p}{r_1 p - 1}}\big)^{(d - 1)\left(1 -  \frac{1}{\theta}\right)} \asymp
$$
\begin{equation}\label{I_2fin2}
\asymp 2^{-m\left(r_1 -\frac{1}{p}\right)}2^{-2m} m^{-2(d - 1)} m^{(d - 1)\left(1 -  \frac{1}{\theta}\right)} \ll 2^{-m r_1}  m^{(d - 1)\left(\frac{1}{p} - \frac{1}{\theta}\right)}.
\end{equation}

Таким чином, згідно з \eqref{I_2os}--\eqref{I_2fin2} отримуємо таке співвідношення
\begin{equation}\label{sum2-fin}
I_2 \ll 2^{-m r_1}  m^{(d - 1)\left(\frac{1}{p} - \frac{1}{\theta}\right)}.
\end{equation}

Для оцінки величини $I_1$ нам знадобиться допоміжне твердження.

Нехай $\mathcal{T}(\Delta Q_n)_q$ позначає одиничну $L_q$-кулю у  просторі поліномів $\mathcal{T}(\Delta Q_n)$. Крім того покладемо
$$
\gamma(q, a, b) =
 \begin{cases}
              \left(\frac{b}{a}\right)^{\frac{1}{q}}\left[ \ln \left(1+\frac{b}{a}\right)\right]^{\frac{1}{q}-\frac{1}{2}}, & a \leq b,  \\
               e^{-\frac{a}{b}}, & a > b.
 \end{cases}
$$

\bf Лема B \rm \cite{Kashin_Temlyakov_1999}\textbf{.}  \it Для $1<q\leq 2$ має місце оцінка
$$
 \varepsilon_M\big(\mathcal{T}(\Delta Q_n)_q, QC\big) \ll n^{\frac{1}{2}} \gamma(q, M, \mathcal{K}|\Delta Q_n|).
$$
\rm

($\mathcal{K}=\mathcal{K}(d)$; інші константи у цій нерівності також не залежать ні від $M$,  ні від $n$).

Отже, представимо величину  $I_1$ у вигляді
\begin{equation}\label{sum1-1}
I_1 = \sum\limits_{n \leq m} \varepsilon_{M_n}\big(S_{\Delta Q_n}(B^{\boldsymbol{r}}_{p,\theta}), QC\big)
+\sum\limits_{m < n \leq m_1} \varepsilon_{M_n}\big(S_{\Delta Q_n}(B^{\boldsymbol{r}}_{p,\theta}), QC\big).
\end{equation}
Далі, згідно з \eqref{delta-fin}  і лемою~В знаходимо
\begin{equation}\label{sum1-1-1}
\sum\limits_{n \leq m} \varepsilon_{M_n}\big(S_{\Delta Q_n}(B^{\boldsymbol{r}}_{p,\theta}), QC\big) \ll \sum\limits_{n \leq m} 2^{-nr_1} n^{(d-1)\left(\frac{1}{p}-\frac{1}{\theta}\right)}n^{\frac{1}{2}}e^{-\frac{M_n}{\mathcal{K}|\Delta Q_n|}}\ll 2^{-mr_1} m^{(d-1)\left(\frac{1}{p}-\frac{1}{\theta}\right)}m^{\frac{1}{2}}.
\end{equation}

Аналогічно при  $n<m\leq m_1$ одержимо
$$
\sum\limits_{m < n \leq m_1} \varepsilon_{M_n}\big(S_{\Delta Q_n}(B^{\boldsymbol{r}}_{p,\theta}), QC\big) \ll
$$
$$
\ll \sum\limits_{m < n \leq m_1} 2^{-nr_1} n^{(d-1)\left(\frac{1}{p}-\frac{1}{\theta}\right)}n^{\frac{1}{2}}\left(\frac{|\Delta Q_n|}{M_n}\right)^{\frac{1}{p}}\left[ \ln \left(1+\frac{|\Delta Q_n|}{M_n}\right)\right]^{\frac{1}{p}-\frac{1}{2}}\ll
$$
\begin{equation}\label{sum1-1-2}
\ll  2^{-mr_1} m^{(d-1)\left(\frac{1}{p}-\frac{1}{\theta}\right)} m^{\frac{1}{2}}.
\end{equation}

Отже, згідно з \eqref{sum1-1}--\eqref{sum1-1-2}  для оцінки величини $I_1$ можемо записати:
\begin{equation}\label{sum1-fin}
I_1 \ll 2^{-mr_1} m^{(d-1)\left(\frac{1}{p}-\frac{1}{\theta}\right)} m^{\frac{1}{2}}.
\end{equation}

Тепер, об'єднавши оцінки \eqref{epsilon-sum}, \eqref{sum2-fin} і  \eqref{sum1-fin} та беручи до уваги, що $M \asymp 2^m m^{d-1}$, для випадку $1<p \leq 2$, $p\leq\theta<\infty$ одержимо
\begin{equation}\label{epsilon-1fin}
\varepsilon_M\big(B^{\boldsymbol{r}}_{p,\theta}, QC \big) \ll M^{-{r_1}}\big(\log^{d - 1} M\big)^{r_1 +\frac{1}{p} - \frac{1}{\theta}} \sqrt{\log M}.
\end{equation}

Далі, скориставшись оцінкою \eqref{epsilon-1fin} одержимо оцінки величини $\varepsilon_M\big(B^{\boldsymbol{r}}_{p,\theta}, QC \big)$ у випадках, що залишилися нерозгянутими.

Нехай $1<p\leq 2$, $1\leq \theta< p$. Тоді  врахувавши, що $B^{\boldsymbol{r}}_{p,\theta}\subset B^{\boldsymbol{r}}_{p,p}$, згідно з \eqref{epsilon-1fin} маємо
\begin{equation}\label{epsilon-2fin}
\varepsilon_M\big(B^{\boldsymbol{r}}_{p,\theta}, QC \big) \ll \varepsilon_M\big(B^{\boldsymbol{r}}_{p,p}, QC \big) \ll M^{-{r_1}}\big(\log^{d - 1} M\big)^{r_1} \sqrt{\log M}.
\end{equation}

Нехай $2<p\leq \infty$, $2 < \theta<\infty$. Тоді $B^{\boldsymbol{r}}_{p,\theta}\subset B^{\boldsymbol{r}}_{2,\theta}$ і тому, скориставшись \eqref{epsilon-1fin}, можемо записати
\begin{equation}\label{epsilon-3fin}
\varepsilon_M\big(B^{\boldsymbol{r}}_{p,\theta}, QC \big) \ll \varepsilon_M\big(B^{\boldsymbol{r}}_{2,\theta}, QC \big) \ll M^{-{r_1}}\big(\log^{d - 1} M\big)^{r_1+\frac{1}{2} - \frac{1}{\theta}} \sqrt{\log M}.
\end{equation}

Нехай $2<p\leq \infty$, $1 \leq \theta\leq 2$. Тоді $B^{\boldsymbol{r}}_{p,\theta}\subset B^{\boldsymbol{r}}_{p,2}\subset B^{\boldsymbol{r}}_{2,2}$ і згідно з \eqref{epsilon-1fin} при $\theta=p=2$ одержимо
\begin{equation}\label{epsilon-4fin}
\varepsilon_M\big(B^{\boldsymbol{r}}_{p,\theta}, QC \big) \ll \varepsilon_M\big(B^{\boldsymbol{r}}_{2,2}, QC \big) \ll M^{-{r_1}}\big(\log^{d - 1} M\big)^{r_1} \sqrt{\log M}.
\end{equation}

Об'єднавши \eqref{epsilon-1fin}--\eqref{epsilon-4fin} приходимо до шуканої оцінки.

Теорему~1 доведено.
\vskip 3 mm

У наступному твердженні встановимо оцінку знизу величини $\varepsilon_M\big(B^{\boldsymbol{r}}_{\infty,\theta}, QC \big)$.

\vskip 3 mm

\bf Теорема 2.  \rm  \it Нехай $r_1 > 0$, $1 \leq \theta < \infty$. Тоді при $d\geq 2$ справедлива оцінка
\begin{equation}\label{Em-QC2}
\varepsilon_M\big(B^{\boldsymbol{r}}_{\infty,\theta}, QC\big) \gg  M^{-{r_1}}(\log^{d - 1} M)^{r_1 + \frac{1}{2} - \frac{1}{\theta}} \sqrt{\log M}.
\end{equation} \rm

{\textbf{\textit{Доведення.}}} Нехай $N_\varepsilon(F,\mathscr{X})$~---   мінімальна  кількість замкнутих куль радіуса ${\varepsilon>0}$  простору $\mathscr{X}$ необхідних для компактного покриття множини $F$,  а $M_\varepsilon(F,\mathscr{X})$~--- максимальна кількість таких точок $x_i \in F$, що $\|x_i - x_j\|_{\mathscr{X}} > \varepsilon$, $i \neq j$.  Тоді справедливі нерівності (див., наприклад, \cite{Kolmogorov_Tixomirov-59})
\begin{equation}\label{Ne-Me}
N_\varepsilon(F,\mathscr{X}) \leq M_{\varepsilon}(F,\mathscr{X}) \leq N_{\frac{\varepsilon}{2}}(F,\mathscr{X}).
\end{equation}

Далі для парних $n$ і  $d\geq 2$ позначимо
$$
Y^d_n = \Big\{\boldsymbol{s}\colon \boldsymbol{s}= (2 l_1, \dots, 2 l_d), l_1 + \dots + l_d = \frac{n}{2}, \boldsymbol{l}\in \mathbb{Z}^d_+\Big\},
$$
$$
\mathscr{D}_n = \bigcup\limits_{\boldsymbol{s} \in Y^d_n}\rho(\boldsymbol{s})
$$
і $\mathcal{T}_r(\mathscr{D}_n)$~---  простір дійсних тригонометричних поліномів $t\in \mathcal{T}(\mathscr{D}_n)$. При цьому зауважимо, що для кількості елементів множин $\mathscr{D}_n $  справедливе співвідношення $|\mathscr{D}_n | \asymp 2^n  n^{d-1}$.

У \cite{Kashin_Temlyakov_1999} (див., також, \cite{Temlyakov_1995JC}) для кожного $n$ побудовано набір функцій $\{f^n_i\}^{A_n}_1$, ${f^n_i \in \mathcal{T}_r(\mathscr{D}_n)}$  з властивостями:
\begin{align} \label{f_i_n}
1) \ \ & \|\delta_{\boldsymbol{s}}(f^n_i)\|_{\infty} \leq 1, & \boldsymbol{s} \in Y^d_n; \notag
\\
2) \ \ &\|f^n_i - f^n_j\|_{QC} \geq C(d) n^{\frac{d}{2}}, & i \neq j;
\\
3) \ \ &  A_n \geq 2^{\frac{|\mathscr{D}_n|}{2}}. \notag
\end{align}

Покажемо, що кожна функція з множини
$$
F_n=\Big\{ C(r,\theta)2^{-nr_1}n^{-\frac{d-1}{\theta}}f^n_i  \Big\}^{A_n}_{i=1}
$$
з деякою константою $C(r,\theta)$ належить класу $B^{\boldsymbol{r}}_{\infty,\theta}$, $1\leq \theta < \infty$.

Маємо
$$
\|f^n_i\|_{B^{\boldsymbol{r}}_{\infty,\theta}} \asymp \sum\limits_{\boldsymbol{s}\in \mathscr{D}_n} 2^{(\boldsymbol{s},\boldsymbol{r})\theta} \|A_{\boldsymbol{s}}(f^n_i)\|^{\theta}_{\infty}=
$$
$$
= \left(\sum\limits_{\boldsymbol{s}\in \mathscr{D}_n} 2^{(\boldsymbol{s},\boldsymbol{r})\theta} \bigg\|A_{\boldsymbol{s}} \ast \sum\limits_{\|\boldsymbol{s}-\boldsymbol{s}'\|\leq 1} \delta_{\boldsymbol{s}'}(f^n_i)\bigg\|^{\theta}_{\infty}\right)^{\frac{1}{\theta}} \leq
$$
$$
\leq \left(\sum\limits_{\boldsymbol{s}\in \mathscr{D}_n} 2^{(\boldsymbol{s},\boldsymbol{r})\theta} \|A_{\boldsymbol{s}}\|^{\theta}_1 \bigg\| \sum\limits_{\|\boldsymbol{s}-\boldsymbol{s}'\|\leq 1} \delta_{\boldsymbol{s}'}(f^n_i)\bigg\|^{\theta}_{\infty}\right)^{\frac{1}{\theta}} \ll
$$
$$
\ll \left(\sum\limits_{\boldsymbol{s}\in \mathscr{D}_n} 2^{(\boldsymbol{s},\boldsymbol{r})\theta} \Bigg( \sum\limits_{\|\boldsymbol{s}-\boldsymbol{s}'\|\leq 1} \|\delta_{\boldsymbol{s}'}(f^n_i)\|_{\infty}\Bigg)^{\theta}\right)^{\frac{1}{\theta}} \ll 2^{nr_1}n^{\frac{d-1}{\theta}}
$$

Отже, $F_n\subset B^{\boldsymbol{r}}_{\infty,\theta}$.

Тепер, беручи до уваги \eqref{Ne-Me} і  скориставшись властивістю 2) для функцій з \eqref{f_i_n} можемо записати
$$
\varepsilon_M\big(B^{\boldsymbol{r}}_{\infty,\theta}, QC \big) \gg \varepsilon_M\big(F_n, QC \big) \gg
$$
$$
\gg 2^{-nr_1}n^{-\frac{d-1}{\theta}}  n^{\frac{d}{2}} = 2^{-nr_1}n^{(d-1)\left(\frac{1}{2}-\frac{1}{\theta}\right)}n^{\frac{1}{2}} \asymp   M^{-{r_1}}(\log^{d - 1} M)^{r_1 + \frac{1}{2} - \frac{1}{\theta}} \sqrt{\log M}
$$

Оцінку \eqref{Em-QC2} встановлено. Теорему~2 доведено.
\vskip 3 mm

З результатів теорем~1 та 2 легко одержати таке твердження.

\vskip 3 mm

\bf Теорема 3.  \rm  \it Нехай $2\leq p \leq \infty$, $2 \leq \theta < \infty$, $r_1 > \frac{1}{2}$. Тоді при $d\geq 2$ справедлива оцінка
\begin{equation}\label{Em-QC3}
\varepsilon_M\big(B^{\boldsymbol{r}}_{p,\theta}, QC\big) \asymp  M^{-{r_1}}(\log^{d - 1} M)^{r_1 + \frac{1}{2} - \frac{1}{\theta}} \sqrt{\log M}.
\end{equation} \rm

{\textbf{\textit{Доведення.}}} Оскільки для $1\leq p <\infty$ має місце вкладення $B^{\boldsymbol{r}}_{\infty,\theta}\subset B^{\boldsymbol{r}}_{p,\theta}$, то скориставшись оцінкою \eqref{Em-QC2}, зокрема, і для $2\leq p \leq \infty$, $2 \leq \theta < \infty$, $r_1 > \frac{1}{2}$,  маємо
\begin{equation}\label{Em-QC3'}
\varepsilon_M\big(B^{\boldsymbol{r}}_{p,\theta}, QC \big) \gg \varepsilon_M\big(B^{\boldsymbol{r}}_{\infty,\theta}, QC \big) \gg M^{-{r_1}}(\log^{d - 1} M)^{r_1 + \frac{1}{2} - \frac{1}{\theta}} \sqrt{\log M}.
\end{equation}

Співставивши  \eqref{Em-QC3'} з оцінкою \eqref{Em-QC} з теореми~1, одержуємо оцінку \eqref{Em-QC3}.

Теорему~3 доведено.

\vskip 5 mm

На завершення  роботи зробимо декілька коментарів.

Одержані у теоремі~3 точні за порядком оцінки ентропійних чисел класів Нікольського--Бєсова $B^{\boldsymbol{r}}_{p,\theta}$ у просторі $QC$ доповнюють відповідні результати  для класів Соболева $W^{\boldsymbol{r}}_{p,\alpha}$ та Нікольського $H^{\boldsymbol{r}}_{p}$, $1 < p \leq \infty$, встановлені Б.\,С.~Кашиним і В.\,М.~Темляковим~\cite{Kashin_Temlyakov_1999}.

Стосовно оцінок ентропійних чисел класів Нікольського--Бєсова $B^{\boldsymbol{r}}_{p,\theta}$ у просторі $L_{\infty}$ зазначимо, що у роботі \cite{Romanyuk_2019AM} встановлено їхні точні за порядком оцінки, але лише у  двовимірному випадку $d=2$. Для зручності порівняння цих оцінок з результатом теореми~3 наведемо відповідне твердження.

\bf Теорема~Б \rm \cite{Romanyuk_2019AM}\textbf{.}    \it Нехай $d=2$, $2\leq p \leq \infty$,  $\boldsymbol{r}=(r_1,r_1)$,  $r_1 > \frac{1}{2}$. Тоді для $2 \leq \theta < \infty$ справедлива оцінка
\begin{equation}\label{Em-Linfty}
\varepsilon_M\big(B^{\boldsymbol{r}}_{p,\theta}, L_{\infty}\big) \asymp  M^{-{r_1}}(\log M)^{r_1 + 1 - \frac{1}{\theta}}.
\end{equation} \rm

Зауважимо, що у випадку $p=\infty$ дана оцінка справедлива при $r_1>0$.

Отже, співставивши \eqref{Em-Linfty} з результатом теореми~3 при  $d=2$ бачимо, що справедливе співвідношення:
$$
\varepsilon_M\big(B^{\boldsymbol{r}}_{p,\theta}, QC \big)\asymp \varepsilon_M\big(B^{\boldsymbol{r}}_{p,\theta}, L_{\infty}\big) \asymp  M^{-{r_1}}(\log M)^{r_1 + 1 - \frac{1}{\theta}},
$$
$$
2\leq p \leq \infty, \ r_1 > \frac{1}{2}, \ 2 \leq \theta < \infty.
$$

Однак питання про порядок величини $\varepsilon_M\big(B^{\boldsymbol{r}}_{p,\theta}, L_{\infty}\big)$, $1\leq \theta \leq \infty$, $1\leq p \leq \infty$ у випадку $d>2$ залишається відкритим.

\vskip 3 mm

\textbf{Contact information:}
Department of the Theory of Functions, Institute of Mathematics of National
Academy of Sciences of Ukraine, 3, Tereshenkivska st., 01024, Kyiv, Ukraine.

\vskip 3 mm

E-mail:  \href{mailto:romanyuk@imath.kiev.ua }{romanyuk@imath.kiev.ua }

\hspace{1.35 cm}
\href{mailto:Yan.Sergiy@gmail.com}{yan.sergiy@gmail.com}

\end{document}